\documentclass[12pt]{article}
\usepackage{amssymb,amsmath,amsthm}
\newcommand{\mr}[1]{\mathrm{#1}}
\newcommand{\mf}[1]{\mathfrak{#1}}
\newcommand{\mc}[1]{\mathcal{#1}}
\newcommand{\eq}[1]{(\ref{#1})}           

\newcommand{\z}{{\bf Z}}                    
\newcommand{\q}{{\bf Q}}
\newcommand{\bq}{\bar{\q}}
\newcommand{\qp}{{\bf Q}_p}
\newcommand{\zp}{{\bf Z}_p}
\newcommand{\ra}{\rightarrow}

\newcommand{\gr}{gr}
\newcommand{\ab}{\mr{ab}}

\DeclareMathOperator{\Gal}{Gal}

\newtheorem{theorem}{Theorem}                            
\newtheorem{prop}[theorem]{Proposition}
\newtheorem{lemma}[theorem]{Lemma}
\newtheorem{cor}[theorem]{Corollary}

\begin{document}
\title{Relationships between conjectures on the
structure of pro-$p$ Galois groups unramified outside $p$}
\author{Romyar T. Sharifi}
\maketitle

\section{Introduction}
The fixed field $\Omega^*$ of the canonical 
representation $\phi \colon G_{\q} \ra \mr{Out}(\pi_1)$ with
\[ 
  \pi_1 = \pi_1^{\mr{pro-}p}({\bf P}^1(\bq) \setminus \{ 0,1,\infty \}) 
\]
is a pro-$p$ extension of $\q(\zeta_p)$ unramified outside $p$
for any prime number $p$ \cite{ihara-ann}.  We study, for odd primes $p$, 
the structure of $G = \Gal(\Omega^*/\q(\zeta_{p^{\infty}}))$ 
and of a certain graded $\zp$-Lie algebra $\mf{g}$ associated to $\phi$
and arising from a filtration of $G$.
To that effect, this article can be viewed as an extension of the article of 
Ihara \cite{ihara-pre} (esp., Lecture I).  

Our primary insight comes from the examination of the relationship between 
elements of $\Gal(\Omega/\q(\zeta_p))$, 
where $\Omega$ is the maximal pro-$p$ extension of $\q(\zeta_p)$ unramified 
outside $p$, and elements of $\Gal(\Omega/\q(\zeta_{p^{\infty}}))$.
More specifically, we construct elements 
$\sigma_m \in \Gal(\Omega/\q(\zeta_{p^{\infty}}))$ 
restricting nontrivially to elements of the $m$th graded pieces 
$\gr^m \mf{g}$ for odd $m \ge 3$.  The elements $\sigma_m$ are obtained 
recursively starting from 
elements of $\Gal(\Omega/\q(\zeta_p))$ which satisfy the property that their 
images in the maximal abelian quotient
generate its odd eigenspaces under the action of $\Gal(\q(\zeta_p)/\q)$.  
In fact, these elements of $\Gal(\Omega/\q(\zeta_p))$ provide suitable 
$\sigma_m$ for odd $m$ with $3 \le m \le p$.
In this way, we are able to employ knowledge of the Galois group of 
$\Omega/\q(\zeta_p)$ in studying the structure of $\mf{g}$.
In particular, it will follow from a consequence of a conjecture of 
Greenberg's in multivariable Iwasawa theory \cite{green-iwa} that $\mf{g}$ 
is not free on the restrictions of the $\sigma_m$ for a large class of 
irregular primes.

Define $\Omega^*_{m-1}$ as the fixed field of the kernel of
\[ \phi_m \colon G_K \ra \mr{Out}(\pi_1/\pi_1(m+1)), \]
where $K = \q(\zeta_{p^{\infty}})$ and 
$\pi_1(m+1)$ denotes the $(m+1)$st term in the lower central series of $\pi_1$.
The graded $\zp$-Lie algebra $\mf{g}$ is defined by setting
$\gr^m \mf{g} = \Gal(\Omega^*_m/\Omega^*_{m-1})$ for $m \ge 1$. 
Recall that, as a $\Gal(K/\q)$-module, 
$\gr^m \mf{g} \cong \zp(m)^{\oplus r_m}$ for some $r_m \ge 0$ 
\cite{ihara-ann,ihara-pre}.  

Consider the filtration $F^m G = \Gal(\Omega^*/\Omega^*_{m-1})$ of $G$ giving
rise to $\mf{g}$.  For each odd integer $m \ge 3$, we let 
$\sigma_m$ denote an element of $F^m G$ with the property that 
$\kappa_m(\sigma_m)$ generates the image of the $m$th Soul\'e character 
$\kappa_m$ on $F^mG$ \cite{soule,ihara-pre}.  In Section \ref{constr}, 
we shall make a particular choice of the
elements $\sigma_m$.  By abuse of notation, the element of $\gr^m\mf{g}$ 
given by the restriction of $\sigma_m$ to $\Omega^*_m$ will also be denoted 
by $\sigma_m$ (and is nontrivial \cite{ihara-msri,ihara-pre}).  We remark 
that $\kappa_m$ and $\kappa_{k}$ induce the same character modulo $p^n$ if
$m \equiv k \bmod p^{n-1}(p-1)$.

Let $S$ denote the free pro-$p$ group on infinitely many generators $s_m$ 
with $m$ odd $\ge 3$, and let $\mf{s}$ denote the free graded $\zp$-Lie 
algebra with generators also denoted $s_m$ in odd degree $\ge 3$.  
Consider the homomorphisms $\Psi \colon S \ra G$ and 
$\psi \colon \mf{s} \ra \mf{g}$ determined by $s_m \mapsto \sigma_m$.
Ihara attributes to Deligne the conjecture that $\psi\otimes \qp$ is an
isomorphism \cite{ihara-msri,ihara-icm,ihara-pre}.  
In fact, Hain and Matsumoto have recently shown 
$\psi \otimes \qp$ to be surjective \cite{hm}.   

We will study the surjectivity, or lack thereof, of $\psi$ itself.
This appears to be a finer question, its answer depending upon arithmetic 
properties of the prime $p$.  We shall show that the validity of
certain conjectures would imply that the map $\psi$ is not surjective 
exactly when $p$ is an irregular prime.  Without any assumptions, our 
results imply that $\psi$ is not an isomorphism for a large class of 
irregular primes.

The following theorem will result from our definition of the elements 
$\sigma_m$, but holds for any choice of $\sigma_m$.

\begin{theorem} \label{reg}
  Let $p$ be an odd regular prime.  Then the homomorphism $\Psi$ is 
  surjective.  If Deligne's conjecture holds for $p$, then $\psi$ and $\Psi$ 
  are isomorphisms and $\Omega = \Omega^*$.
\end{theorem}

Additionally, we shall obtain a refinement of this in the form of an 
extension of Theorem I-2(ii) of \cite{ihara-pre}.
It will also follow from our methods that if $p$ is regular and 
$\Omega = \Omega^*$, then the map $\Psi$ is injective.  Furthermore, we shall 
show that injectivity 
of $\Psi$ together with surjectivity of $\psi$ forces $\psi$ to be 
injective.  Hence, we obtain the following result.

\begin{theorem} \label{delreg}
    Let $p$ be an odd regular prime.  If $\Omega = \Omega^*$ and $\psi$ is
    surjective, then Deligne's conjecture holds for $p$.
\end{theorem}

For any number field $F$, Greenberg has a conjecture regarding the structure 
of the Galois group of a large pro-$p$ extension of $F$ unramified outside 
primes above $p$ \cite{green-iwa} (see also \cite{green-gal}).  Let 
$F_{\infty}$ 
denote the compositum of all $\zp$-extensions of $F$ (which is necessarily 
unramified outside $p$), and let $L_{\infty}$ denote the maximal abelian 
unramified pro-$p$ extension of $F_{\infty}$.  Greenberg's conjecture 
states that $X = \Gal(L_{\infty}/F_{\infty})$ has annihilator of height at 
least $2$ as a module over the multivariable Iwasawa algebra 
$\Lambda = \zp[[ \Gal(F_{\infty}/F) ]]$ \cite{green-iwa}.
By class field theory, $X$ is seen to be isomorphic to the inverse limit 
$A_{\infty}$ of the $p$-parts of the ideal class groups of finite 
subextensions of $F$ in $F_{\infty}$.

We consider Greenberg's conjecture for the field $F = \q(\zeta_p)$, which we 
refer to as Greenberg's conjecture for the prime $p$.  Since $A_{\infty} = 0$ 
if $p$ is regular, Greenberg's conjecture holds trivially for regular primes.  
Let $M_{\infty}$ denote the maximal abelian pro-$p$ extension of $F_{\infty}$ 
unramified outside $p$.  Greenberg's conjecture for $p$ is equivalent to 
$\Gal(M_{\infty}/F_{\infty})$ having no torsion as a module over $\Lambda$
\cite{mc,ln}.
In particular, Greenberg's conjecture for an irregular prime $p$ implies 
that $\Gal(\Omega/F)$ has no free pro-$p$ quotient on $(p+1)/2$ generators 
\cite{mc} (see also \cite{ln}).
The absence of such a free pro-$p$ quotient is the key to the following 
theorem.

\begin{theorem} \label{main}
    Let $p$ be an irregular prime.  If Greenberg's conjecture holds for 
    $p$, then $\psi$ and $\Psi$ are not isomorphisms.  In particular, if 
    Deligne's conjecture also holds for $p$, then $\psi$ and $\Psi$ are not 
    surjective.
\end{theorem}

McCallum has proven Greenberg's conjecture for a large class of irregular 
primes \cite{mc}.  These are exactly those primes $p$ for which both 
the $p$-part $A$ of the ideal class group of 
$\q(\zeta_p)$ and $(U/\bar{E})[p^{\infty}]$ are of order $p$, where
$U$ denotes the unit group of $\q_p(\zeta_p)$ and $\bar{E}$ denotes the
closure of the image of $\z[\zeta_p]^*$ in $U$
(see \cite{mc,marshall} for equivalent conditions).  David Marshall 
has extended this result to a class of primes $p$ for which $A$ which may be 
cyclic of any positive $p$-power order \cite{marshall}.  Currently, there 
are no known examples in which $A$ is cyclic of order greater than $p$.  
However, there is an abundance of primes for which $A$ is not cyclic, for 
instance $p = 157$ and $691$.

\section{Construction of the elements} \label{constr}

Choose any element $\tau \in \Gal(\Omega/\q)$ that restricts to a
generator of $\Gal(K/\q)$.  Let $\delta = \lim_{i \ra \infty} 
\tau^{p^i}$, an element of order $p-1$ restricting to a generator of 
$\Gal(F/\q)$, and let $\gamma = \tau^{p-1}$, which 
commutes with $\delta$ and restricts to a generator of $\Gal(K/F) \cong \zp$.  
By abuse of notation, we 
also denote by $\delta$ and $\gamma$ the restrictions to subfields 
Galois over $\q$.  
For any $m \in \z$, let 
$\epsilon_m \in \zp[\Gal(\Omega/\q)]$ be 
the element
\[ 
  \epsilon_m = \frac{1}{p-1} \sum_{i=0}^{p-2} \chi(\delta^i)^{-m} 
  \delta^i, 
\]
where $\chi \colon G_{\q} \ra \zp^*$ denotes the cyclotomic 
character.  

Let $E$ denote a pro-$p$ extension of $F$ (unramified outside $p$)
which is also Galois over $\q$. 
If $E/F$ is nonabelian, we cannot in general define an action of the
idempotent $\epsilon_m$ on $\Gal(E/F)$.  However, the following provides 
something of a substitute.  For $g \in \Gal(E/F)$, we define
\begin{equation} \label{emact} 
   g^{\epsilon_m} = (g \cdot \delta g^{\chi(\delta)^{-m}} 
   \delta^{-1} \cdot \delta^2 g^{\chi(\delta^2)^{-m}} \delta^{-2} 
  \cdots \delta^{p-2} g^{\chi(\delta^{p-2})^{-m}} 
  \delta^{-p+2})^{1/(p-1)}.
\end{equation}
Note that $g^{\epsilon_m} = g^{\epsilon_{m'}}$ whenever 
$m \equiv m' \bmod p-1$.
We also define $g^{\epsilon_m^i}$ to be the $i$th iterate
$(\ldots(g^{\epsilon_m})^{\epsilon_m}\ldots)^{\epsilon_m}$.  
Although $\epsilon_m^2 = \epsilon_m$, we do not necessarily have
$g^{\epsilon_m^2} = g^{\epsilon_m}$.
Instead, we have the following lemma.

\begin{lemma} \label{gm}
    For any $g \in \Gal(E/F)$, the element
    \[ g^{(m)} = \lim_{i \ra \infty} g^{\epsilon_m^{i}} \]
    is well-defined and satisfies
    \begin{equation} \label{deltagm}
        \delta g^{(m)} \delta^{-1} = (g^{(m)})^{\chi(\delta)^{m}}.
    \end{equation}
\end{lemma}

\begin{proof}
    Set $N = \Gal(E/F)$.  Let $x \in N$, and
    let $B$ denote the normal closure in $\Gal(E/\q)$ of the 
    subgroup generated by the elements $[x,\delta^j x \delta^{-j}]$
    with $1 \le j \le p-2$.
    We begin by proving the claim that
    \[ \delta x^{\epsilon_m} \delta^{-1} 
    (x^{\epsilon_m})^{-\chi(\delta)^m} \in B. \]
    To see this, note that by definition \eq{emact} we have
    \begin{equation} \label{claim} 
       \delta x^{\epsilon_m} \delta^{-1}
       = (\delta x \delta^{-1} \cdot \delta^2 x^{\chi(\delta)^{-m}} 
       \delta^{-2} \cdots \delta^{p-2} x^{\chi(\delta^{p-3})^{-m}} 
       \delta^{-p+2} \cdot x^{\chi(\delta)^{m}})^{1/(p-1)}.
    \end{equation}
    In $N/B$, the terms of \eq{claim} commute, and the right-hand side of 
    equation \eq{claim} equals 
    \[
       (x^{\chi(\delta)^m})^{\epsilon_m} = 
       (x^{\epsilon_m})^{\chi(\delta)^m},
    \]
    which proves the claim.
    
    Let $N(i)$ denote the $i$th term in the lower 
    central series of $N$.  We now show inductively that for $i \ge 0$ and
    $1 \le j \le p-2$ we have
    \begin{equation} \label{a_i}
      a_{i,j} = \delta^j g^{\epsilon_m^i} \delta^{-j} 
      (g^{\epsilon_m^i})^{-\chi(\delta^j)^{m}}
      \in N(i+1).
    \end{equation}
    Note that $a_{0,j} \in N = N(1)$, and assume that
    $a_{i-1,j} \in N(i)$ for each $j$ and some $i \ge 1$.  Our earlier 
    claim implies that
    the element $a_{i,j}$ is contained in the normal closure $C(i+1)$ in 
    $\Gal(E/\q)$ of the group generated by the commutators 
    \[
       b_{i,l} = [g^{\epsilon_m^{i-1}},\delta^l 
       g^{\epsilon_m^{i-1}} \delta^{-l}] = 
       [g^{\epsilon_m^{i-1}},a_{i-1,l}]
    \]
    with $1 \le l \le p-2$.
    Hence $b_{i,l} \in N(i+1)$ for each $l$ and therefore $a_{i,j} \in N(i+1)$.
    
    Similarly, we have that
    \begin{equation} \label{c_i}
      c_i = g^{\epsilon_m^i}(g^{\epsilon_m^{i-1}})^{-1} \in C(i+1) 
      \triangleleft N(i+1).
    \end{equation}
    As $\cap_{i=1}^{\infty} N(i) = \{1\}$, we have that
    $g^{(m)}$ is well-defined by \eq{c_i}.  Furthermore, \eq{deltagm} holds
    by \eq{a_i} with $j=1$.
\end{proof}

We shall now make our choice of the elements $\sigma_m$ described in 
the introduction.  First,  for odd $m$ 
with $3 \le m \le p$, we choose any element $t_m \in \Gal(\Omega^*/K)$ 
such that $\kappa_m(t_m)$ generates $\kappa_m(\Gal(\Omega^*/K))$, and we set 
$g_m = t_m^{(m)}$, which also has maximal image.  For such $m$, set 
$\sigma_m = g_m$.  We recursively define the other elements $g_m$ and 
$\sigma_m$ by
\begin{equation} \label{sm}
    g_{m+p-1} = \gamma g_m \gamma^{-1} g_m^{-\chi(\gamma)^m}
\end{equation}
and
\begin{equation} \label{sigmam} 
   \sigma_{m+p-1} = (\gamma \sigma_m \gamma^{-1} 
   \sigma_m^{-\chi(\gamma)^m})^{(m)} 
\end{equation}
for any odd $m \ge 3$.

(We remark that it is not necessary to take the limit element in 
\eq{sigmam}; a finite iteration 
\[ 
  (\gamma \sigma_m \gamma^{-1} 
  \sigma_m^{-\chi(\gamma)^m})^{\epsilon_m^{i}} 
\]
would work as an element $\sigma_{m+p-1}$, in fact with $i \le 1$ if 
$m \ge p-2$.)

Let $L^*$ denote the maximal abelian extension of $K$ in $\Omega^*$.

\begin{lemma} \label{restab}
    The elements $\sigma_m$ and $g_m$ have the same restriction to
    $L^*$ for all odd $m \ge 3$.
\end{lemma}

\begin{proof}
    For $3 \le m \le p$, the statement is trivially true.  
    Note that $\epsilon_m$ defines an idempotent endomorphism of 
    $G^{\ab} = \Gal(L^*/K)$.  Also, $\epsilon_m$ commutes with 
    conjugation by $\gamma$ on $G^{\ab}$.  Hence we have 
    \[ 
      \sigma_{m+p-1}|_{L^*} =  (\gamma \sigma_m\gamma^{-1}  
      \sigma_m^{-\chi(\gamma)^m})^{\epsilon_m}|_{L^*}
      = \gamma \sigma_m^{\epsilon_m} \gamma^{-1}  
      (\sigma_m^{\epsilon_m})^{-\chi(\gamma)^m}|_{L^*}
      = \gamma \sigma_m \gamma^{-1} \sigma_m^{-\chi(\gamma)^m}|_{L^*}. 
    \]
    The statement now follows immediately by induction.
\end{proof}

Let $v_p$ denote the $p$-adic valuation on $\zp$.

\begin{prop} \label{kapsig}
    Let $m$ denote an odd integer with $m \ge 3$.
    The element $\sigma_m$ fixes $\Omega_{m-1}^*$, and its image under
    $\kappa_m$ generates $\kappa_m(F^mG)$.  Furthermore, if $m = k + j(p-1)$ 
    for some integers $k \ge 3$ and $j \ge 0$, then
    \[ 
      v_p(\kappa_{m}(\sigma_{m})) = v_p((jp)!) +
      v_p(\kappa_{m}(\sigma_k)). 
    \]
\end{prop}

\begin{proof}
    We first show that $\sigma_m$ fixes $\Omega_{m-1}^*$ inductively.
    Set $\Omega^*_{k} = K$ for $k < 0$.  Assume that $\sigma_m$ fixes
    $\Omega_{m-p+1}^*$.  Once we show that $\sigma_m$ fixes $\Omega_{m-1}^*$, 
    then $\gamma\sigma_{m}\gamma^{-1}\sigma_{m}^{-\chi(\gamma)^{m}}$ 
    will fix $\Omega_m^*$, as $\gr^m \mf{g}$ has Tate twist $m$.
   
    We prove the claim that if $x$ fixes $\Omega^*_{m-p+1}$ 
    then $y = x^{(m)}$ fixes $\Omega^*_{m-1}$.  Note first that $y$ 
    fixes $\Omega_{m-p+1}^*$, so we assume inductively that it fixes 
    $\Omega_{k-1}^*$ with $m-p+2 \le k \le m-1$.  As $\Omega^*_k$ is abelian 
    over $\Omega^*_{k-1}$ and $\mf{g}_k$ has Tate twist $k$, we see that 
    \[ 
      \delta y \delta^{-1}|_{\Omega_k^*} = y^{\chi(\delta)^k}|_{\Omega_k^*}. 
    \]
    From Lemma \ref{gm}, we therefore have that $y$ 
    will fix $\Omega^*_k$, since $k \not\equiv m \bmod p-1$.
    We conclude that $y$ fixes $\Omega^*_{m-1}$.  The first statement of
    the proposition now follows.

    Note that
    \begin{equation} \label{kapconj}
        \kappa_m(\gamma t \gamma^{-1}) = \chi(\gamma)^m\kappa_m(t)
    \end{equation}
    for any $t \in G$.  Therefore, we see that
    \begin{equation} \label{compval}
       v_p(\kappa_m(\gamma t \gamma^{-1} t^{-\chi(\gamma)^l})) 
       = v_p((m-l)p) + v_p(\kappa_m(t))
    \end{equation}
    for any $t \in G$.
    The last statement of the proposition now follows
    from the recursive definition \eq{sm} of $g_m$ and \eq{compval}.        

    By Lemma \ref{restab}, the element $g_m$ has the
    same image under $\kappa_m$ as $\sigma_m$. 
    Let $\tilde{L}$ denote the abelian subextension of $\Omega$ generated 
    by roots of cyclotomic $p$-units, and
    note that $\kappa_m$ can be considered as a homomorphism of
    $A = \Gal(\tilde{L}/K)$ \cite{ihara-pre}.  
    To prove the second statement of the proposition, it
    suffices to show that $\kappa_m(g_m)$ generates 
    $\kappa_m(\Gal(\tilde{L}/\tilde{L} \cap \Omega^*_{m-1}))$.
    
    Let $i$ denote the unique integer with $3 \le i \le p$ and 
    $m \equiv i \bmod p-1$.  Let $A_i = A^{\epsilon_i}$, 
    $B = \Gal(\tilde{L}/F)$, and $B_i = (B^{\ab})^{\epsilon_i}$. 
    Let $h_i$ denote an element of $A_i$ restricting to a generator of 
    the procyclic group $B_i$ \cite[Ch.\ 8,13]{wash}.  Then
    $h_i$ topologically generates $A_i$ as a normal subgroup of $B$.
    
    Let $\tilde{x}$ denote the restriction of an element 
    $x \in \Gal(\Omega/F)$ to $\tilde{L}$. 
    We claim that $A_i$ is also the normal closure in $B$ of the procyclic
    subgroup generated by $\tilde{g_i}$. If not, then since 
    $\tilde{g_i} \in A_i$, we must have
    \[ \tilde{g_i} = x^p[\gamma,y] \]
    for some $x,y \in A_i$.  Clearly, we would then have
    \[ v_p(\kappa_i(g_i)) > v_p(\kappa_i(h_i)), \]
    contradicting the definition of $g_i$.
 
    Let $A_m$ denote the largest normal subgroup of 
    $A_i$ fixing $\Omega^*_{m-1} \cap \tilde{L}$.  We have $\kappa_m(F^mG) 
    = \kappa_m(A_m)$, and so we need only show that $A_m$ is the normal 
    closure in $B$ of the procyclic subgroup generated by 
    $\tilde{g}_m$.  Inductively assuming this is true for $k$, we prove it 
    for $k+p-1$.  If $x \in A_k$, then since $A$ is abelian, 
    $x \in A_k$ is a product of conjugates of $\tilde{g}_k$ by powers of 
    $\tilde{\gamma}$.  From equation \eq{kapconj}, we see that
    \[ \kappa_k(\prod \gamma^j g_k^{a_j} \gamma^{-j}) = 0 \]
    if and only if
    \[ \sum {a_j \chi(\gamma)^{kj}} = 0. \]
    Observing \eq{sm}, we see that any such element is a product of 
    conjugates of $\tilde{g}_{k+p-1}$ by powers of $\tilde\gamma$.
\end{proof}

Proposition \ref{kapsig} has an interesting application to the study
of the stable derivation algebra $\mc{D}$ over $\z$ considered by Ihara 
\cite{ihara-groth,ihara-isr,ihara-pre}.   
Ihara has shown that there is a canonical embedding of graded 
$\zp$-Lie algebras
\[ \iota \colon \mf{g} \hookrightarrow \mc{D} \otimes \zp \]
and has conjectured that $\gr^m \iota$ is an isomorphism for (at least)
$m < p$.  There is also a canonical map
\[ \lambda_m \colon \gr^m \mc{D} \ra \z \]
(denoted $\gr^m(c)$ in \cite{ihara-pre})
which, after extending $\zp$-linearly, we may compose with 
$\gr^m\iota$  to obtain a map $\lambda_m^{(p)}$.
The latter map is related to the Soul\'e character $\kappa_m$ by the
formula \cite{ihara-msri,ihara-pre}
\begin{equation} \label{relate}
    \kappa_m = (p^{m-1}-1)(m-1)!\lambda_m^{(p)}
\end{equation}
on $F^mG$.
Let $N_m$ denote the positive generator of the image ideal 
of $\lambda_m$.

\begin{cor}
    Let $p$ denote an odd prime satisfying Vandiver's conjecture.
    Let $m \ge 3$ be an odd positive integer, and let $k$
    denote the largest integer less than or equal to $(m-3)/(p-1)$.
    Then 
    \begin{equation} \label{nm}
        v_p(N_m) \le v_p((kp)!) - v_p((m-1)!),
    \end{equation} 
    with equality if $\gr^m \iota$ is an isomorphism.    
\end{cor}

\begin{proof}
    Since $\lambda_m^{(p)}(F^mG) \subseteq N_m\zp$, this follows directly 
    from \eq{relate}, the surjectivity of $\kappa_{m-k(p-1)}$ under
    Vandiver's conjecture \cite{is} and the last 
    statement of Proposition \ref{kapsig}.
\end{proof}

\section{Proofs of the main results}

Before proving any theorems regarding $\psi$ and $\Psi$, we make the 
following general points.

\begin{lemma} \label{freegp}
Fix $r \ge 1$, and consider a pro-$p$ group $\mc{F}$ topologically generated by 
elements $y$ and $x_i$ with $1 \le i \le r$.  
For each $1 \le i \le r$, define $x_{i,1} = x_i$ and 
\begin{equation} \label{free}
   x_{i,j+1} = y x_{i,j} y^{-1} x_{i,j}^{-1+pa_{i,j}} = [y,x_{i,j}] 
   x_{i,j}^{pa_{i,j}}  
\end{equation}
for some $a_{i,j} \in \zp$ for each $j \ge 1$.
Let $H$ denote the normal closure of the pro-$p$ subgroup of $\mc{F}$ generated
by the $x_i$ with $1 \le i \le r$.  Then:
\begin{enumerate}
    \item[a.] The elements $x_{i,j}$ with $1 \le i \le r$ and $j \ge 1$ 
    topologically generate $H$.
    \item[b.] If $\mc{F}$ is a free pro-$p$ group on the elements $y$ 
    and $x_i$ with $1 \le i \le r$, then $H$ is a free 
    pro-$p$ group on the elements $x_{i,j}$.
    \item[c.] Assume $H$ is a free pro-$p$ group on the 
    elements $x_{i,j}$ with $1 \le i \le r$ and $j \ge 1$. If
    $y^{p^n} \notin H$ for all $n$, then $\mc{F}$ is a free pro-$p$ group 
    on the elements $y$ and $x_i$ with $1 \le i \le r$.
\end{enumerate}
\end{lemma}

\begin{proof}
    To prove part (a), we need only show that, for each $k \ge 0$, 
    the element $y^k x_i y^{-k}$ is contained in the pro-$p$ subgroup 
    generated by the elements $x_{i,j}$.  This follows easily from the fact 
    that the group generated by $x_{i,j+1}$ and $x_{i,j}$ contains 
    $y x_{i,j} y^{-1}$.
    
    For any pro-$p$ group $N$, let $N[j]$ denote the $j$th term in its 
    descending central $p$-series, and set $\bar{N} = N/N[2]$.
    Now assume $\mc{F}$ is free pro-$p$ on $y$ and the $x_i$.
    Note that $H$ is a free pro-$p$ group as a closed subgroup of a free 
    pro-$p$ group. 
    Hence it is free on the $x_{i,j}$ if and only if the images 
    of the $x_{i,j}$ form a minimal generating set of $\bar{H}$.
    
    Fix $i$, and let $D = D_i$ denote the free pro-$p$ subgroup of $\mc{F}$ 
    generated by $x_i$ and $y$.  Let $C = C_i$ denote the normal closure in 
    $D$ of the pro-$p$ subgroup generated by $x_i$.  By freeness of 
    $\mc{F}$, we have
    \[ \bar{H} \cong \bigoplus_{k=1}^{r} \bar{C}_k. \]
    Hence we are reduced to showing that the images of the $x_{i,j}$ form a 
    minimal generating set of $\bar{C}$.  
    As $\bar{C}$ is elementary abelian, we have an injection
    \[ 
       \bar{C} \hookrightarrow \bigoplus_{j \ge 0} 
       (C \cap D[j])/(C \cap D[j+1]). 
    \]
    Finally, as $x_{i,j} \in D[j] \setminus D[j+1]$, we obtain minimality,
    proving part (b).
    
    As for part (c), consider a free presentation of $\mc{F}$ on generators 
    $\tilde{x}_i$ and $\tilde{y}$ mapping to $x_i$ and $y$:
    \[ 1 \ra R \ra \tilde{\mc{F}} \ra \mc{F} \ra 1. \]
    By parts (a) and (b), the group $\tilde{\mc{F}}$ has a free 
    subgroup $\tilde{H}$ on elements $\tilde{x}_{i,j}$ defined as in
    \eq{free}.  Since $H$ is free on the $x_{i,j}$, the map $\tilde{H} \ra H$ 
    is an isomorphism.  Hence $R \cap \tilde{H} = 0$.  Thus $R$ is isomorphic
    to $R\tilde{H}/\tilde{H}$, which is a subgroup of the group
    $\tilde{\mc{F}}/\tilde{H} \cong \zp$ generated by the 
    image of $\tilde{y}$.  Hence, if $R$ is 
    nontrivial, then $\tilde{y}^{p^n} \in R\tilde{H}$ for some $n \ge 0$.  
    Therefore we have $y^{p^n} \in H$, a contradiction.  Hence $R$ 
    is trivial, proving 
    (c).
\end{proof}

We now prove Theorem \ref{reg}.

\begin{proof}[Proof of Theorem \ref{reg}]
    For $p$ regular, it is well-known that $\Gal(\Omega/F)$ is free pro-$p$ 
    on the generators $\gamma$ and $g_m$ with $3 \le m \le p$ and 
    $m$ odd.  By the construction \eq{sm} of the $g_m$ and parts (a) and (b) 
    of Lemma \ref{freegp}, it 
    follows that $H = \Gal(\Omega/K)$ is freely generated as 
    a pro-$p$ group by the $g_m$ with $m$ odd $\ge 3$.  
    As $H$ is pronilpotent, and since by Lemma \ref{restab} the
    elements $\sigma_m$ and $g_m$ 
    have the same image on the Galois group of the 
    maximal abelian subextension of $K$ in $\Omega$, we have that
    $H$ is also free on the $\sigma_m$.  As $G = \Gal(\Omega^*/K)$ 
    is a quotient of $H$, the map $\Psi$ is surjective.  (Note that
    the surjectivity of $\Psi$ does not depend on the choice of 
    elements $\sigma_m$.)
    
    Assuming Deligne's conjecture, we have that $\psi$, and hence 
    $\Psi$, is injective.  Therefore, $\Psi$ is an isomorphism and
    $\Omega = \Omega^*$.  Finally, the injectivity of $\psi$ plus 
    the surjectivity of $\Psi$ yield the surjectivity of $\psi$ 
    \cite[\S III-I-2]{ihara-pre}.
\end{proof}

The following serves as a rough converse to the injectivity implies 
surjectivity theorem for $\psi \otimes \qp$ proven by Ihara 
\cite[Theorem I-1]{ihara-pre}.  Recall that $\mf{s}$ arises as the
graded Lie algebra associated to a filtration $F^m S$ on $S$ compatible 
with the filtration on $G$ \cite{ihara-pre}.

\begin{theorem} \label{surinj}
If $\Psi$ is injective and $\psi$ is surjective, then $\psi$ is 
injective.
\end{theorem}

\begin{proof}
    We begin by assuming merely that $\Psi$ is injective and $\psi$ is not.  
    Then $\psi(\bar{x}) = 0$ for some 
    nonzero $\bar{x} \in \gr^m\mf{s}$ and $m \ge 3$.  Choose a lift 
    $x_1 \in S$ of $\bar{x}$.  Since $\Psi$ is injective, the image element 
    $g_1 = 
    \Psi(x_1)$ is nontrivial.  Let $k_1$ be maximal such that $g_1$ fixes 
    $\Omega^*_{k_1-1}$.  Let $\bar{g}_1$ be the image of $g_1$ in 
    $gr^{k_1} \mf{g}$.  
    By \cite{hm}, we have that $\psi \otimes \qp$ is surjective, so there 
    exists  $\bar{y}_1 \in 
    gr^{k_1} \mf{s}$ such that $\psi(\bar{y}_1) = 
    p^{n_1} \bar{g}_1$ for some minimum possible $n_1 \ge 0$.  Set $x_2 
    = x_1^{p^{n_1}} y_1^{-1}$ for $y_1$ lifting 
    $\bar{y}_1$.  By construction and the injectivity of $\Psi$, there exists
    $k_2 > k_1$ maximal such that $g_2 = \Psi(x_2)$ fixes $\Omega^*_{k_2-1}$.
    By induction, we obtain a sequence of elements $x_i \in S$, each of which 
    restricts to some multiple of $\bar{x} \in \gr^m \mf{s}$, and corresponding 
    sequences of elements $y_i$ and exponents $n_i$ such that
    $x_{i+1} = x_i^{p^{n_i}}y_i^{-1}$.
    
    If $\psi$ is surjective, then each $n_i$ is zero.  Since the 
    sequence of numbers $k_i$ is increasing,
    the sequence $x_i$ has a limit $x \in S$ which
    restricts to $\bar{x} \in \gr^m \mf{s}$.  We see that 
    $\Psi(x) = 0$, as $\Psi(x)$ will fix $\Omega^*_k$ for every $k$.  
    This is a contradiction, proving the theorem.  
    Note that, removing the assumption of the surjectivity of $\psi$, this 
    argument yields that infinitely many of the $n_i$ are non-zero.
\end{proof}

Theorem \ref{delreg} now follows as a corollary to Theorem \ref{surinj}.

\begin{proof}[Proof of Theorem \ref{delreg}]
    If $\Omega = \Omega^*$, then since $p$ is regular, Lemma \ref{freegp}b 
    implies that $G$ is free on the generators $g_m$ with $m$ odd $\ge 3$
    and hence is free on the $\sigma_m$ by Lemma \ref{restab}.
    That is, $\Psi$ is an isomorphism.  Since by assumption $\psi$ is
    surjective, Theorem \ref{surinj} implies that $\psi$ is an 
    isomorphism.
\end{proof}

Before proving Theorem \ref{main}, we make the following remark.

\begin{lemma} \label{Psiso}
  The state of $\Psi$ being an isomorphism does not depend on the choice
  of elements $\sigma_m$.  
\end{lemma}

\begin{proof}
  For each odd $m \ge 3$, let $\tau_m$ be an element of $F^mG$ with maximal 
  possible image under $\kappa_m$.  Let $\Psi' \colon S \ra G$ be a homomorphism 
  satisfying $\Psi'(s_m) = \tau_m$ for each odd $m \ge 3$. 
  Assume that $\Psi$ is an isomorphism.  The surjectivity of $\Psi$
  forces that
  \[ 
    \tau_m|_{L^*} = \prod_{\substack{k \ge m \\ k \mr{\ odd}}} 
    \sigma_k^{a_{m,k}}|_{L^*}
  \]
  for some $a_{m,m} \in \zp^*$ and $a_{m,k} \in \zp$ for odd $k > m$.
  Setting $a_{m,k} = 0$ for $k < m$, we find that the matrix formed
  by the $a_{m,k}$ with $m$ and $k$ odd $\ge 3$ is upper triangular 
  and invertible.  Then the $\tau_m|_{L^*}$ form a minimal generating set of
  $G^{\ab} = \Gal(L^*/K)$, and hence the $\tau_m$ freely 
  generate $G$, as $G$ is a free pro-$p$ group.
\end{proof}

We may now prove Theorem \ref{main}.

\begin{proof}[Proof of Theorem \ref{main}]
    The idea of the proof is to show that if $\psi$ or $\Psi$ is an 
    isomorphism then $\mc{G} = \Gal(\Omega^*/F)$ is free, 
    contradicting the corollary of Greenberg's conjecture that
    $\Gal(\Omega/F)$ has no free pro-$p$ quotient of rank $(p+1)/2$.
    The state of $\psi$ or $\Psi$ being an isomorphism does not depend on the
    choice of generators $\sigma_m$ by \cite[\S III-6]{ihara-pre} and Lemma 
    \ref{Psiso}, so we use our 
    previously defined generators from \eq{sigmam} (which is possible by 
    Proposition \ref{kapsig}).  
    
    If $\psi$ is an isomorphism, then $\Psi$ is as well 
    \cite[\S III-6]{ihara-pre}.  Assume that $\Psi$ is an isomorphism.
    This means that the 
    $\sigma_m$ freely generate $G$, so Lemma \ref{restab} implies that the 
    $g_m$ do as well.  We clearly have that $\mc{G}$ is generated by $\gamma$
    and $g_m$ with $m$ odd and $3 \le m \le p$.  As $\gamma^{p^n} \notin G$ 
    for any $n$,
    it follows by Lemma \ref{freegp}c that $\mc{G}$ is
    freely generated by these elements, finishing the proof.
\end{proof}

\section{The filtration on the Galois group}

We now consider the question of where in the filtration of $G$ the 
nonsurjectivity of $\Psi$ may first occur.
In this regard, we also have the following extension of Theorem I-2(ii) of 
\cite{ihara-pre} which removes the assumption $m < p$. 
It can also be viewed as a more precise version of the contrapositive to 
Theorem \ref{reg} that if $\Psi$ is not surjective, then $p$ is irregular. 

Let $G_m = \Gal(\Omega^*_{m}/K)$ and set $S_m = S/F^{m+1}S$ \cite{ihara-pre}.
We also recall the map $\Psi_m \colon S_m \ra G_m$ with 
$\Psi_m(s_i) = \sigma_i$ for $3 \le i \le m$ and $i$ odd.

\begin{theorem} \label{ihgen}
  Let $p$ be an odd prime and assume that
  $\Psi$ is not surjective.  If $m$ is minimal such that $\Psi_m$ is not
  surjective, then $\Omega^*_m$ contains a nontrivial 
  elementary abelian extension of $F$ with Tate twist $m$
  which is linearly disjoint from the fixed field of the kernel of 
 $\kappa_m$ if $m$ is odd.  If $m$ is even, then $p$ divides $B_m$ and 
  $L^* \cap \Omega^*_m$ has a nontrivial even part.  If $m$ is odd, then
  $p$ divides $B_{p-m}$ and Vandiver's conjecture fails at $p$.
\end{theorem}

\begin{proof}
  We remark that the surjectivity, or lack thereof, of $\Psi$ and $\Psi_m$ 
  is independent of the choice of the elements $\sigma_k$.
  Thus, Proposition \ref{kapsig} again allows us to use our choice of 
  $\sigma_k$ from \eq{sigmam}.
 
  Set $B = \gr^m \mf{g}$ if $m$ is even, and let $B$ denote the kernel of
  $\kappa_m$ on $\gr^m \mf{g}$ if $m$ is odd.  Set $A = B/(B \cap G_m(2))$.
  Similarly to the 
  proof of Theorem I-2(ii) of \cite{ihara-pre}, we have an isomorphism of 
  abelian groups
  \[ G_m^{\ab} \cong \Psi(S_m)^{\ab} \times A. \]

  We claim that $\Psi(S_m)^{\ab}$ is a 
  $\zp^*$-submodule of $G_m^{\ab}$.  (Ihara points out that a similar group 
  need not be $\zp^*$-stable.) In fact, we have by definition that
  \[  
    \delta \sigma_k \delta^{-1} = \sigma_k^{\chi(\delta)^k}
  \]
  and 
  \[  
    \gamma \sigma_k \gamma^{-1}|_{L^*} =
    \sigma_k^{\chi(\gamma)^i}\sigma_{k+p-1}|_{L^*}.
  \] 
  This proves the claim, so $A$ is a direct summand of $G_m^{\ab}$ as a 
  $\zp^*$-submodule.  

  Therefore, $A$ gives rise to an abelian pro-$p$ 
  extension $\Sigma$ of $K$ unramified outside $p$ with Tate twist $m$ 
  and satisfying $\Sigma \cap \tilde{L} = K$.  
  The maximal elementary abelian subextension of $\Sigma/K$ descends to the 
  desired elementary abelian $p$-extension of $F$.
\end{proof}

One may ask if the assumption $m < p$ in Ihara's result actually needs to be 
removed.  That is, is the smallest $m$ such that $\Psi_m$ is not surjective 
also the smallest $m$ such that $p$ divides $B_m$?  We shall briefly describe
what an answer should involve.

Let us assume Vandiver's conjecture at $p$.  That $p$ divides $B_m$ indicates 
the existence of a relation in $\mc{G}$.  If $m$ is any number such that 
$\Psi_{m-1}$ is surjective, then the relation induced in $G_m$ can be put in 
the form
\[ h^{p^c} = \Psi_m(s) \]
for some $h \in G_m$, $c \ge 1$, and $s \in [S_m,S_m]$
which is not a $p$th power of a nontrivial element.
The map $\Psi_m$ is not surjective if and only if $\Psi_m(s) \neq 1$.  
The surjectivity of $\Psi_m$ is therefore governed by both the structure 
of the relation in $\mc{G}$ and the form of the filtration on $G$, two 
objects about which we require more information.  

For example, in the case $p = 691$ and $m = 12$ considered by Ihara 
\cite{ihara-pre} we know that the commutators $[\sigma_3,\sigma_9]$ and 
$[\sigma_5,\sigma_7]$ are linearly independent in
$\gr^{12} \mf{g} \otimes \q_{691}$  \cite{matsumoto}, and so we have 
information about the filtration.  However, to prove the nonsurjectivity 
of $\Psi_{12}$ one still needs to demonstrate that the relation in $\mc{G}$
actually ``involves'' these commuators so that in $G_{12}$ it reduces 
to a nontrivial relation of the form
\[ h^{691} = [\sigma_3,\sigma_9]^a[\sigma_5,\sigma_7]^b. \]

\end{document}